\documentclass[11pt]{amsart}

\usepackage{etex}
\usepackage{amsmath, amssymb}
\usepackage{array}
\usepackage[frame,cmtip,arrow,matrix,line,graph,curve]{xy}
\usepackage{graphpap, color, paralist, pstricks}
\usepackage[mathscr]{eucal}
\usepackage[pdftex]{graphicx}
\usepackage[pdftex,colorlinks,backref=page,citecolor=blue]{hyperref}
\usepackage{tikz}

\newtheorem{thm}{Theorem}[section]
\newtheorem{prop}[thm]{Proposition}
\newtheorem{cor}[thm]{Corollary}
\newtheorem{lem}[thm]{Lemma}

\theoremstyle{plain}
\newtheorem{defn}[thm]{Definition}

\newtheorem{rem}[thm]{Remark}

\newcommand{\cc}{\mathbb{C}}
\newcommand{\pp}{\mathbb{P}}

\newcommand{\Gr}{\mathrm{Gr}}
\newcommand{\PGL}{\mathrm{PGL}}
\newcommand{\End}{\mathrm{End}}

\newcommand{\EEnd}{\mathscr E\!nd}

\newcommand{\Hilb}{\mathrm{Hilb}}

\newcommand{\GH}{\mathrm{GH}}
\newcommand{\GPB}{\mathrm{GPB}}
\newcommand{\GPH}{\mathrm{GPH}}
\newcommand{\GGPH}{\mathrm{GGPH}}
\newcommand{\GV}{\mathrm{GV}}
\newcommand{\id}{\mathrm{id}}

\newcommand{\Spec}{\mathrm{Spec}\,}
\newcommand{\bSpec}{\textbf{\textrm{Spec}} \,}
\newcommand{\md}{\mathrm{mod}}
\newcommand{\Sym}{\mathrm{Sym}}
\newcommand{\rk}{\mathrm{rank} \,}
\newcommand{\cE}{\mathcal{E} }
\newcommand{\cF}{\mathcal{F} }
\newcommand{\cG}{\mathcal{G}}
\newcommand{\cM}{\mathcal{M} }
\newcommand{\cO}{\mathcal{O} }
\newcommand{\cU}{\mathcal{U} }
\newcommand{\cV}{\mathcal{V} }
\newcommand{\cW}{\mathcal{W} }
\newcommand{\cY}{\mathcal{Y} }
\newcommand{\fS}{\mathfrak{S} }
\newcommand{\tG}{\tilde{G}}
\newcommand{\s}{\sigma }
\newcommand{\om}{\omega }
\newcommand{\dss}{\displaystyle}

\begin{document}

\title[A geometry of moduli of Higgs pairs on an irreducible nodal curve of arithmetic genus one]{A geometry of the moduli space of Higgs pairs on an irreducible nodal curve of arithmetic genus one}

\author{Sang-Bum Yoo}
\address{Department of Mathematics Education, Gongju National University of Education, 27 Ungjin-ro, Gongju-si, Chungcheongnam-do, 32553, Republic of Korea}
\email{sbyoo@gjue.ac.kr}


\begin{abstract}	
We describe the moduli space of Higgs pairs on an irreducible nodal curve of arithmetic genus one and its geometric structures in terms of the Hitchin map and a flat degeneration of the moduli space of Higgs bundles on an elliptic curve.
\end{abstract}

\maketitle

\section{Introduction}

\subsection{Motivations and results}

The moduli space of Higgs bundles on a smooth curve has been intensively studied in view of a mirror symmetry that was first raised by T. Hausel and M. Thaddeus in \cite{HT03}. Our work in this paper provides concrete ingredients to be useful later when we prove or disprove a mirror symmetry phenomenon over a singular curve.

The purpose of this paper is to describe the moduli space of Higgs pairs on an irreducible nodal curve of arithmetic genus one and its geometric structures in terms of the Hitchin map and a flat degeneration of the moduli space of Higgs bundles on an elliptic curve explicitly.

Throughout this paper, $Y$ denotes a reduced irreducible projective curve of arithmetic genus one, with only one ordinary node $p$, defined over $\cc$, let $\nu:X\to Y$ be the normalization map and let $\nu^{-1}(p)=\{p_{1},p_{2}\}$. Note that $X\cong\pp^{1}$.

There are two ways to study Higgs pairs on $Y$ like torsion-free sheaves on $Y$ (See \cite{Gie84,Se82}). One is to compare them with generalized parabolic Higgs bundles on $X$ in order to give the Hitchin map on the moduli space of Higgs pairs on $Y$ (See \cite{Bh14}). Another is to compare them with Gieseker-Hitchin pairs on $X$ attached with a chain of projective lines in order to give a flat degeneration of the moduli space of Higgs bundles on an elliptic curve (See \cite{BBN13,BBNarx13}).

For a positive integer $n$ and an arbitrary integer $d$, we consider the following moduli spaces.
\begin{itemize}
\item Let $U_{Y}(n,d)$ be the moduli space of semistable torsion-free sheaves of rank $n$ and degree $d$ on $Y$.
\item Let $U_{X}^{\GPB}(n,d)$ be the moduli space of semistable generalized parabolic vector bundles of rank $n$ and degree $d$ on $X$.
\item Let $\cM_{Y}(n,d)$ be the moduli space of semistable Higgs pairs of rank $n$ and degree $d$ on $Y$.
\item Let $\cM_{X}^{\GGPH}(n,d)$ be the moduli space of semistable good generalized parabolic Higgs bundles of rank $n$ and degree $d$ on $X$.
\end{itemize}

We describe $\cM_{Y}(n,d)$ and $\cM_{Y}^{\GGPH}(n,d)$ as follows.

\begin{thm}[Theorem \ref{coprime moduli of Higgs pairs}]\label{thm1}
If $\gcd(n,d)=1$, then $\cM_{Y}(n,d)\cong U_{Y}(n,d)\times\cc\cong Y\times\cc$.
\end{thm}

\begin{thm}[Theorem \ref{noncoprime Higgs pairs}]\label{thm2}
If $\gcd(n,d)>1$, then there is no stable Higgs pairs of rank $n$ and degree $d$ on $Y$.
\end{thm}

\begin{thm}[Theorem \ref{noncoprime moduli of Higgs pairs}]\label{thm3}
Let $\gcd(n,d)=h$.
\begin{enumerate}
\item There exists a bijective morphism
$$\Sym^{h}(Y\times\cc)\to\cM_{Y}(n,d).$$
\item $\cM_{Y}(n,d)$ is irreducible.
\end{enumerate}
\end{thm}

Theorem \ref{thm1} follows from a simplified stability of Higgs pairs on $Y$ and a description of $U_{Y}(n,d)$ by \cite{BB08} in the case $\gcd(n,d)=1$. The proof of Theorem \ref{thm2} is given by considering a degeneration of the moduli space of semistable Higgs bundles of rank $n$ and degree $d$ on an elliptic curve and using the nonexistence of stable Higgs bundles of rank $n$ and degree $d$ on an elliptic curve in the case $\gcd(n,d)>1$. Theorem \ref{thm3} follows from a family of semistable Higgs pairs of rank $n$ and degree $d$ on $Y$ parametrized by $(Y\times\cc)\times\overset{h}{\cdots}\times(Y\times\cc)$ induced from the universal family of stable Higgs pairs of rank $\dss\frac{n}{h}$ and degree $\dss\frac{d}{h}$ on $Y$ parametrized by $Y\times\cc$ in the case $\gcd(n,d)=h$.

\begin{thm}[Theorem \ref{coprime moduli and stable locus of noncoprime moduli of GPHs}]\label{thm4}
\begin{enumerate}
\item\label{thm4-1} If $\gcd(n,d)=1$, then $\cM_{X}^{\GGPH}(n,d)\cong U_{X}^{\GPB}(n,d)\times\cc\cong\pp^{1}\times\cc$.
\item\label{thm4-2} If $\gcd(n,d)>1$, then the stable locus $\cM_{X}^{\GGPH}(n,d)^{s}$ of $\cM_{X}^{\GGPH}(n,d)$ is empty.
\end{enumerate}
\end{thm}

\begin{thm}[Theorem \ref{noncoprime moduli of good GPHs}]\label{thm5}
Let $\gcd(n,d)=h$.
\begin{enumerate}
\item There exists a bijective morphism
$$\Sym^{h}(\pp^{1}\times\cc)\to\cM_{X}^{\GGPH}(n,d).$$
\item $\Sym^{h}(\pp^{1}\times\cc)$ is the normalization of $\cM_{X}^{\GGPH}(n,d)$. 
\item $\cM_{X}^{\GGPH}(n,d)$ is irreducible.
\end{enumerate}
\end{thm}

Theorem \ref{thm4}-(\ref{thm4-1}) follows from a simplified stability of generalized parabolic Higgs bundles on $X$ and the known result that $U_{X}^{\GPB}(n,d)$ is the normalization of $U_{Y}(n,d)$ in the case $\gcd(n,d)=1$. Theorem \ref{thm4}-(\ref{thm4-2}) follows from the nonexistence of stable Higgs bundles of rank $n$ and degree $d$ on an elliptic curve in the case $\gcd(n,d)>1$. Theorem \ref{thm5} follows from the normality of $\Sym^{h}(\pp^{1}\times\cc)$ and a family of semistable generalized parabolic Higgs bundles of rank $n$ and degree $d$ on $X$ parametrized by $(\pp^{1}\times\cc)\times\overset{h}{\cdots}\times(\pp^{1}\times\cc)$ induced from the universal family of stable generalized parabolic Higgs bundles of rank $\dss\frac{n}{h}$ and degree $\dss\frac{d}{h}$ on $X$ parametrized by $\pp^{1}\times\cc$ in the case $\gcd(n,d)=h$.

We also describe all fibers of the Hitchin maps $H$ on $\cM_{Y}(n,d)$ and $H^{\GGPH}$ on $\cM_{X}^{\GGPH}(n,d)$ as follows.

\begin{thm}[Corollary \ref{fiber of h of GGPH}]\label{thm6}
Let $\gcd(n,d)=h$. The generic fiber of the Hitchin map $H^{\GGPH}$ on $\cM_{X}^{\GGPH}(n,d)$ is set-theoretically isomorphic to $\pp^{1}\times\overset{h}{\cdots}\times\pp^{1}$. The fiber over an arbitrary point of the base is set-theoretically isomorphic to $\pp^{m_{1}}\times\overset{h}{\cdots}\times\pp^{m_{l}}$ where $h=m_{1}+\cdots+m_{l}$. The fiber over an arbitrary point of the base is isomorphic to $\pp^{1}$ for the case $\gcd(n,d)=1$.
\end{thm}

\begin{thm}[Corollary \ref{fiber of h of Higgs pair}]\label{thm7}
Let $\gcd(n,d)=h$. The generic fiber of the Hitchin map $H$ on $\cM_{Y}(n,d)$ is set-theoretically isomorphic to $Y\times\overset{h}{\cdots}\times Y$. The fiber over an arbitrary point of the base is set-theoretically isomorphic to $\Sym^{m_{1}}Y\times\overset{l}{\cdots}\times\Sym^{m_{l}}Y$ where $h=m_{1}+\cdots+m_{l}$. The fiber over an arbitrary point of the base is isomorphic to $Y$ for the case $\gcd(n,d)=1$.
\end{thm}

We define the Hitchin map $H^{\GGPH}$ on $\cM_{X}^{\GGPH}(n,d)$ and then it induces the Hitchin map $H$ on $\cM_{Y}(n,d)$ by a surjective birational morphism $f:\cM_{X}^{\GGPH}(n,d)\to\cM_{Y}(n,d)$. $H^{\GGPH}$ is set-theoretically identified with the projection
$$\pi_{h}:\Sym^{h}(\pp^{1}\times\cc)\to\Sym^{h}\cc,$$
$$[(x_{1},t_{1}),\cdots,(x_{h},t_{h})]_{\fS_{h}}\mapsto[t_{1},\cdots,t_{h}]_{\fS_{h}}.$$
Then the identification of fibers of $H^{\GGPH}$ and $\pi_{h}$ gives Theorem \ref{thm6}. Theorem \ref{thm7} follows from Theorem \ref{thm6} and the surjective birational morphism $f:\cM_{X}^{\GGPH}(n,d)\to\cM_{Y}(n,d)$.

We finally describe a flat degeneration of the moduli space of stable Higgs bundles of rank $n$ and degree $d$ on an elliptic curve for the case $\gcd(n,d)=1$. Let $Z\to T$ be a flat family of irreducible complex projective curves of arithmetic genus one parametrized by a smooth curve $T$ such that for all $t\ne t_{0}$, the fiber $Z_{t}$ is an elliptic curve and $Z_{t_0}\cong Y$. Let $\cG_{Z/T}(n,d)\to T$ be a flat family of the moduli spaces of stable Gieseker-Hitchin pairs of rank $n$ and degree $d$ over $Z$ parametrized by $T$, which is constructed in \cite[Theorem 1]{BBN13} and \cite[Proposition 5.13]{BBNarx13}. Let $H^{\GH}$ be the Hitchin map on $\cG_{Z/T}(n,d)$.

\begin{thm}[Theorem \ref{flat degeneration}]\label{thm8}
\begin{enumerate}
\item $\cG_{Z/T}(n,d)\cong Z\times\cc$ as $T$-schemes.
\item $\nu_{*}^{\GH}:\cG_{Z/T}(n,d)\to\cM_{Z/T}(n,d)$ is identified with the identity map $\id_{Z\times\cc}:Z\times\cc\to Z\times\cc$, where $\cM_{Z/T}(n,d)\to T$ is a flat family of the moduli spaces of stable Higgs pairs of rank $n$ and degree $d$ over $Z$ parametrized by $T$.
\item The fiber of the Hitchin map $H^{\GH}$ on $\cG_{Z/T}(n,d)$ is isomorphic to $Z$.
\end{enumerate}
\end{thm}

\subsection{Organization of the paper}

In Section \ref{tf sheaves and GPBs}, we describe $U_{Y}(n,d)$ and $U_{X}^{\GPB}(n,d)$ explicitly. In Section \ref{Higgs pairs and GHPs}, we describe $\cM_{Y}(n,d)$ and $\cM_{X}^{\GGPH}(n,d)$ explicitly. In Section \ref{Hitchin map}, we give descriptions of all fibers of the Hitchin maps on $\cM_{X}^{\GGPH}(n,d)$ and $\cM_{Y}(n,d)$ respectively. In Section \ref{degeneration}, we prove that $\cG_{Z/T}(n,d)\cong Z\times\cc$ as $T$-schemes and the fiber of the Hitchin map on $\cG_{Z/T}(n,d)$ is isomorphic to $Z$.

\subsection*{Acknowledgements}
The author thanks Young-Hoon Kiem for his comments and suggestions on an earlier draft.

\section{torsion-free sheaves and generalized parabolic vector bundles}\label{tf sheaves and GPBs}

In this section we aim to describe the moduli space of torsion-free sheaves on $Y$ and that of generalized parabolic vector bundles on $X$ explicitly by relating these moduli spaces.

\begin{defn}[\cite{Bh92}]
A {\bf generalized parabolic vector bundle (GPB)} of rank $n$ and degree $d$ on $X$ is a pair $(E,F_{1}(E))$ where $E$ is a vector bundle of rank $n$ and degree $d$ on $X$ and $F_{1}(E)$ is a $n$-dimensional subspace of $E_{p_{1}}\oplus E_{p_{2}}$.
\end{defn}

\begin{defn}[\cite{Bh92}]
A GPB $(E,F_{1}(E))$ is {\bf semistable (}respectively, {\bf stable)} if for every proper subbundle $N\subset E$,
$$\frac{\deg N+\dim(F_{1}(E)\cap(N_{p_{1}}\oplus N_{p_{2}}))}{\rk N}\le(<)\frac{\deg E+\dim F_{1}(E)}{\rk E}.$$
\end{defn}

Consider a GPB $(E,F_{1}(E))$ of rank $n$ and degree $d$ on $X$. To $(E,F_{1}(E))$, we associate a torsion-free sheaf $F$ of rank $n$ and degree $d$ on $Y$ by the following short exact sequence :
$$0\to F\to\nu_{*}E\to\nu_{*}(E)\otimes\cc(p)/F_{1}(E)\to 0.$$

\begin{prop}[Proposition 4.2 of \cite{Bh92}]\label{equivalence of stabilities between GPBs and sheaves}
$(E,F_{1}(E))$ is semistable (respectively, stable) of rank $n$ and degree $d$ if and only if $F$ is semistable (respectively, stable) of rank $n$ and degree $d$.
\end{prop}

Let $U_{Y}(n,d)$ be the moduli space of semistable torsion-free sheaves of rank $n$ and degree $d$ on $Y$ and let $U_{X}^{\GPB}(n,d)$ be the moduli space of generalized parabolic vector bundles of rank $n$ and degree $d$ on $X$. Denote the stable loci of $U_{Y}(n,d)$ and $U_{X}^{\GPB}(n,d)$ by $U_{Y}(n,d)^{s}$ and $U_{X}^{\GPB}(n,d)^{s}$ respectively.

We begin with referring to the results of a classification of stable torsion-free sheaves on $Y$ in \cite{BB08}.

\begin{prop}[Lemma 2.2 and Theorem 2.5 of \cite{BB08}]\label{coprime moduli and stable locus of noncoprim moduli}
\begin{enumerate}
\item\label{coprime moduli and stable locus of noncoprim moduli1} If $\gcd(n,d)=1$, then $U_{Y}(n,d)\cong Y$.
\item\label{coprime moduli and stable locus of noncoprim moduli2} If $\gcd(n,d)>1$, then there is no stable torsion-free sheaf of rank $n$ and degree $d$ over $Y$, that is, $U_{Y}(n,d)^{s}$ is empty.
\end{enumerate}
\end{prop}

Let us denote $\xymatrix{Y\ar[r]^{\cong\quad\quad}&U_{Y}(n,d)}$ of Proposition \ref{coprime moduli and stable locus of noncoprim moduli}-(\ref{coprime moduli and stable locus of noncoprim moduli1}) by $\zeta_{n,d}$.

Next we classify stable GPBs on $X$. Combining Proposition \ref{equivalence of stabilities between GPBs and sheaves} with Proposition \ref{coprime moduli and stable locus of noncoprim moduli}-(\ref{coprime moduli and stable locus of noncoprim moduli2}), we have the following statement.

\begin{lem}\label{noncoprime GPBs}
If $\gcd(n,d)>1$, then there is no stable GPB of rank $n$ and degree $d$ over $X$.
\end{lem}

The following result relates $U_{Y}(n,d)$ to $U_{X}^{\GPB}(n,d)$.

\begin{prop}[Theorem 1 and Theorem 3 of \cite{Bh92}, Proposition 2.1 of \cite{Sun00}]\label{GPB is a normalization}
$U_{X}^{\GPB}(n,d)$ is the normalization of $U_{Y}(n,d)$.
\end{prop}

The following classification is immediately obtained from Proposition \ref{coprime moduli and stable locus of noncoprim moduli}-(\ref{coprime moduli and stable locus of noncoprim moduli1}), Lemma \ref{noncoprime GPBs} and Proposition \ref{GPB is a normalization}.

\begin{prop}\label{coprime moduli and stable locus of noncoprim moduli of GPBs}
\begin{enumerate}
\item\label{coprime moduli and stable locus of noncoprim moduli of GPBs1} If $\gcd(n,d)=1$, then $U_{X}^{\GPB}(n,d)\cong\pp^{1}$.
\item\label{coprime moduli and stable locus of noncoprim moduli of GPBs2} If $\gcd(n,d)>1$, then $U_{X}^{\GPB}(n,d)^{s}$ is empty.
\end{enumerate}
\end{prop}

Denote $\xymatrix{\pp^{1}\ar[r]^{\cong\quad\quad}&U_{X}^{GPB}(n,d)}$ of Proposition \ref{coprime moduli and stable locus of noncoprim moduli of GPBs}-(\ref{coprime moduli and stable locus of noncoprim moduli of GPBs1}) by $\zeta_{n,d}^{\GPB}$. 

Now we classify semistable torsion-free sheaves on $Y$ and semistable GPBs on $X$. By Proposition \ref{coprime moduli and stable locus of noncoprim moduli}-(\ref{coprime moduli and stable locus of noncoprim moduli2}) and Lemma \ref{noncoprime GPBs}, the Jordan-H\"{o}lder filtrations for torsion-free sheaves and GPBs imply the following observations.

\begin{lem}\label{splitting of ss sheaves and GPBs}
Assume $\gcd(n,d)=h$.
\begin{enumerate}
\item\label{splitting of ss sheaves} Any semistable torsion-free sheaf $F$ of rank $n$ and degree $d$ over $Y$ is S-equivalent to $F_{1}\oplus\cdots\oplus F_{h}$, where each $F_{i}$ is stable of rank $\dss\frac{n}{h}$ and degree $\dss\frac{d}{h}$.
\item\label{splitting of ss GPBs} Any semistable GPB $(E,F_{1}(E))$ of rank $n$ and degree $d$ over $X$ is S-equivalent to $(E_{1},F_{1}(E_{1}))\oplus\cdots\oplus(E_{h},F_{1}(E_{h}))$, where each $(E_{i},F_{1}(E_{i}))$ is stable of rank $\dss\frac{n}{h}$ and degree $\dss\frac{d}{h}$.
\end{enumerate}
\end{lem}

When $\gcd(n,d)=1$, there exist universal families of stable torsion-free sheaves and stable GPBs of rank $n$ and degree $d$.

\begin{lem}\label{univ family of sheaves}
Assume that $\gcd(n,d)=1$. Then there exists a universal family $\cV_{n,d}$ of stable torsion-free sheaves of rank $n$ and degree $d$ parametrized by $Y$ such that for every $y\in Y$,
$$\zeta_{n,d}(y)=[(\cV_{n,d})_{y}]_{S},$$
where $[(\cV_{n,d})_{y}]_{S}$ is the S-equivalence class of $(\cV_{n,d})_{y}$.
\end{lem}
\begin{proof}
Since $\chi(E(m))=nm+d$ for any $E\in U_{Y}(n,d)$ and $\gcd(n,d)=1$, by \cite[Corollary 4.6.6]{HL97} and Proposition \ref{coprime moduli and stable locus of noncoprim moduli}-(\ref{coprime moduli and stable locus of noncoprim moduli1}), we prove the statement.
\end{proof}

\begin{lem}\label{univ family of GPBs}
Assume that $\gcd(n,d)=1$. Then there exists a universal family $(\cW_{n,d},F_{1}(\cW_{n,d}))$ of stable GPBs of rank $n$ and degree $d$ parametrized by $X$ such that for every $x\in X$,
$$\zeta_{n,d}^{\GPB}(x)=[(\cW_{n,d},F_{1}(\cW_{n,d}))_{x}]_{S},$$
where $[(\cW_{n,d},F_{1}(\cW_{n,d}))_{x}]_{S}$ is the S-equivalence class of $(\cW_{n,d},F_{1}(\cW_{n,d}))_{x}$.
\end{lem}
\begin{proof}
By Proposition 3.16 of \cite{Bh92} and Proposition \ref{coprime moduli and stable locus of noncoprim moduli of GPBs}-(\ref{coprime moduli and stable locus of noncoprim moduli of GPBs1}), we prove the statement.
\end{proof}

If $\gcd(n,d)=h>1$, $\dss n'=\frac{n}{h}$ and $\dss d'=\frac{d}{h}$, then we can consider the families
$$\cV_{n,d}=\cV_{n',d'}\underset{Y}{\times}\overset{h}{\cdots}\underset{Y}{\times}\cV_{n',d'}$$
of polystable torsion-free sheaves parametrized by $Y\times\overset{h}{\cdots}\times Y$ and
$$(\cW_{n,d},F_{1}(\cW_{n,d}))=(\cW_{n',d'},F_{1}(\cW_{n',d'}))\underset{X}{\times}\overset{h}{\cdots}\underset{X}{\times}(\cW_{n',d'},F_{1}(\cW_{n',d'}))$$
of polystable GPBs parametrized by $\pp^{1}\times\overset{h}{\cdots}\times\pp^{1}$.

The following maps induced by $\cV_{n,d}$ and $(\cW_{n,d},F_{1}(\cW_{n,d}))$,
$$\nu_{\cV_{n,d}}:Y\times\overset{h}{\cdots}\times Y\to U_{Y}(n,d)$$
and
$$\nu_{(\cW_{n,d},F_{1}(\cW_{n,d}))}:\pp^{1}\times\overset{h}{\cdots}\times\pp^{1}\to U_{X}^{\GPB}(n,d),$$
are surjective by Lemma \ref{splitting of ss sheaves and GPBs} and factor through $\Sym^{h}Y$ and $\Sym^{h}\pp^{1}$. Now we complete the classification as follows.

\begin{prop}\label{coprime and noncoprime moduli of sheaves}
\begin{enumerate}
\item For $\gcd(n,d)=h$, there exists a bijective morphism
$$\Sym^{h}Y\to U_{Y}(n,d).$$
\item For $\gcd(n,d)=h$, $U_{X}^{\GPB}(n,d)\cong\Sym^{h}\pp^{1}\cong\pp^{h}$.
\end{enumerate}
\end{prop}
\begin{proof}
\begin{enumerate}
\item $\nu_{\cV_{n,d}}:Y\times\overset{h}{\cdots}\times Y\to U_{Y}(n,d)$ induces a bijective morphism
$$\Sym^{h}Y\to U_{Y}(n,d).$$
\item $\nu_{(\cW_{n,d},F_{1}(\cW_{n,d}))}:\pp^{1}\times\overset{h}{\cdots}\times\pp^{1}\to U_{X}^{\GPB}(n,d)$ induces a bijective morphism
$$\zeta_{n,d}^{\GPB}:\Sym^{h}\pp^{1}\to U_{X}^{\GPB}(n,d).$$
Since $\Sym^{h}\pp^{1}\cong\pp^{h}$ is connected and $U_{X}^{\GPB}(n,d)$ is normal, $\zeta_{n,d}^{\GPB}$ is an isomorphism by Zariski's main theorem.
\end{enumerate}
\end{proof}

\section{Higgs pairs and generalized parabolic Higgs bundles}\label{Higgs pairs and GHPs}
In this section we describe the moduli space of Higgs pairs on $Y$ and that of generalized parabolic Higgs bundles on $X$ explicitly. Note that the dualizing sheaf $\om_{Y}$ is trivial.

\begin{defn}[\cite{Bh14}]
A {\bf Higgs pair} of rank $n$ and degree $d$ on $Y$ is a pair $(E,\phi_{E})$ where $E$ is a torsion-free sheaf of rank $n$ and degree $d$ on $Y$ and $\phi_{E}$ is a global section of $\EEnd E$.
\end{defn}

\begin{defn}[\cite{Bh14}]
A Higgs pair $(E,\phi_{E})$ is {\bf semistable (}respectively, {\bf stable)} if for every proper $\phi_{E}$-invariant subsheaf $N\subset E$,
$$\mu(N)\le(<)\mu(E),$$
where $\dss \mu(E)=\frac{\deg E}{\rk E}$ is the slope of $E$.
\end{defn}

\begin{defn}[\cite{Bh14}]
A {\bf generalized parabolic Higgs bundle (GPH)} of rank $n$ and degree $d$ on $X$ is a triple $(E,\phi_{E},F_{1}(E))$ where $E$ is a vector bundle of rank $n$ and degree $d$ on $X$, $\phi_{E}$ is a global section of $\EEnd E$ and $F_{1}(E)$ is a $n$-dimensional subspace of $E_{p_{1}}\oplus E_{p_{2}}$. A GPH $(E,\phi_{E},F_{1}(E))$ is {\bf good} if $\phi_{E}|_{p_{1}+p_{2}}(F_{1}(E))\subset F_{1}(E)$.
\end{defn}

\begin{defn}[\cite{Bh14}]
A GPH $(E,\phi_{E},F_{1}(E))$ is {\bf semistable (}respectively, {\bf stable)} if for every proper $\phi_{E}$-invariant subbundle $N\subset E$,
$$\frac{\deg N+\dim(F_{1}(E)\cap(N_{p_{1}}\oplus N_{p_{2}}))}{\rk N}\le(<)\frac{\deg E+\dim F_{1}(E)}{\rk E}.$$
\end{defn}

Let $\cM_{Y}(n,d)$ be the moduli space of semistable Higgs pairs of rank $n$ and degree $d$ on $Y$ and $\cM_{Y}(n,d)^{s}$ denotes the stable locus. Let $\cM_{X}^{\GPH}(n,d)$ be the moduli space of semistable GPHs of rank $n$ and degree $d$ on $X$ and $\cM_{X}^{\GPH}(n,d)^{s}$ denotes the stable locus. Let $\cM_{X}^{\GGPH}(n,d)$ be the moduli space of semistable good GPHs of rank $n$ and degree $d$ on $X$ and $\cM_{X}^{\GGPH}(n,d)^{s}$ denotes the stable locus. Note that $\cM_{X}^{\GGPH}(n,d)$ is a closed subscheme of $\cM_{X}^{\GPH}(n,d)$ by Theorem 4.8 of \cite{Bh14}.

The triviality of $\om_{Y}$ gives a simplicity of the study of the semistability of Higgs pairs on $Y$.

\begin{lem}\label{stability}
\begin{enumerate}
\item A Higgs pair $(E,\phi)$ is semistable if and only if $E$ is semistable.
\item If $\gcd(n,d)=1$ and $(E,\phi)\in\cM_{Y}(n,d)$, then $(E,\phi)$ is stable if and only if $E$ is stable.
\end{enumerate}
\end{lem}
\begin{proof}
Since the dualizing sheaf $\om_{Y}$ is trivial, the proof is same as that of Proposition 4.1 of \cite{FGN14}.
\end{proof}

We start with classifying stable Higgs pairs on $Y$.

\begin{thm}\label{coprime moduli of Higgs pairs}
If $\gcd(n,d)=1$, then $\cM_{Y}(n,d)\cong U_{Y}(n,d)\times\cc\cong Y\times\cc$.
\end{thm}
\begin{proof}
Since $\cM_{Y}(n,d)\to\cM_{Y}(n,d+n\deg M),\,(E,\phi_{E})\mapsto(E\otimes M,\phi_{E}\otimes\id_{M})$ is an isomorphism for some fixed line bundle $M$ on $Y$, we may assume that $d>n$ as \cite[page 281]{Nit91}.

Let $U$ be the restriction of the universal sheaf to $Y\times R^{s}$ where $R$ is the open subset of the quot scheme parametrizing quotient sheaves of rank $n$ and degree $d$ in \cite[Section 3]{Nit91}. Applying \cite[Lemma 3.5]{Nit91} with $\cF=\EEnd\,U$, we get a linear scheme $F\to R^{s}$ given by
$$F=\bSpec\Sym_{\cO_{R^{s}}}(\pi_{R^{s}*}\EEnd\,U)^{\vee},$$
where $\pi_{R^{s}}:Y\times R^{s}\to R^{s}$ is the projection onto $R^{s}$. Since $\pi_{R^{s}*}\EEnd\,U\cong\cO_{R^{s}}$ by \cite[Lemma 4.6.3]{HL97}, we have $F\cong R^{s}\times\cc$. It follows from Lemma \ref{stability} that $F^{s}\cong R^{s}\times\cc$.

Hence by the construction of \cite[Section 3]{Nit91} and Proposition \ref{coprime moduli and stable locus of noncoprim moduli}-(\ref{coprime moduli and stable locus of noncoprim moduli1}),
$$\cM_{Y}(n,d)\cong F^{s}/\!/\PGL(d)\cong(R^{s}/\!/\PGL(d))\times\cc\cong U_{Y}(n,d)\times\cc\cong Y\times\cc.$$
\end{proof}

Let us denote $\xymatrix{Y\times\cc\ar[r]^{\cong}&\cM_{Y}(n,d)}$ of Theorem \ref{coprime moduli of Higgs pairs} by $\eta_{n,d}$.

\begin{thm}\label{noncoprime Higgs pairs}
If $\gcd(n,d)>1$, then there is no stable Higgs pairs of rank $n$ and degree $d$ over $Y$. In other words $\cM_{Y}(n,d)^{s}$ is empty.
\end{thm}
\begin{proof}
Let $Z\to S$ be a flat family of irreducible complex projective curves of
arithmetic genus one parametrized by a smooth curve $S$, and let $s_{0}\in S$ be a base point such
that for all $s\ne s_{0}$, the fiber $Z_{s}$ is an elliptic curve and $Z_{s_0}\cong Y$. Let $M\to S$ be the relative moduli variety over $S$ such that $M_{s}$ is the moduli space of semistable Higgs bundles of rank $n$ and degree $d$ on $Z_{s}$ for all $s\ne s_{0}$ and $M_{s_{0}}=\cM_{Y}(n,d)$ (See \cite{Simp94} for the existence of $M$).

Let $M^{s}\subset M$ denote the subset corresponding to stable Higgs pairs. By the openness of the stability condition, $M^s$ is an open subset of $M$ (See Proposition 3.1 of \cite{Nit91}). Assume that there is a stable Higgs pair on $Y$ of rank $n$ and degree $d$ in an irreducible component $M_{\text{irr}}$ of $M$. Then $M^{s}\cap M_{\text{irr}}\cap M_{t_{0}}$ is nonempty, hence $M^{s}\cap M_{\text{irr}}$ is nonempty. Then $M^{s}\cap(M_{\text{irr}}\setminus M_{t_0})$ is a nonempty open subset of $M_{\text{irr}}$. Consequently, $M_s$ is nonempty for some $s\ne s_{0}$. This contradicts the fact that there are no stable Higgs bundles of rank $n$ and degree $d$ on an elliptic curve (\cite[Proposition 4.3]{FGN14} and FACT of \cite{Tu93}). This completes the proof.
\end{proof}

Next we classify stable good GPHs on $X$. We first recall an equivalence of semistabilities between good GPHs on $X$ and Higgs pairs on $Y$. Consider a good GPH $(E,\phi_{E},F_{1}(E))$ of rank $n$ and degree $d$ on $X$. To $(E,\phi_{E},F_{1}(E))$, we associate a Higgs pair $(E,\phi_{E})$ of rank $n$ and degree $d$ on $Y$ by the following commutative diagram of short exact sequences :
\begin{equation}\label{GGPH to Higgs pair}\xymatrix{0\ar[r]&F\ar[r]\ar[d]_{\phi_{F}}&\nu_{*}E\ar[r]\ar[d]_{\nu_{*}\phi_{E}}&\nu_{*}E\otimes\cc(p)/F_{1}(E)\ar[r]\ar[d]_{(\nu_{*}\phi_{E})_{p}}&0\\
0\ar[r]&F\ar[r]&\nu_{*}E\ar[r]&\nu_{*}E\otimes\cc(p)/F_{1}(E)\ar[r]&0.}\end{equation}

\begin{prop}[Theorem 2.9 of \cite{Bh14}]\label{equivalence of semistabilities between good GPHs and Higgs pairs}
A good GPH $(E,\phi_{E},F_{1}(E))$ is semistable if and only if $(F,\phi_{F})$ is semistable.
\end{prop}

Indeed we also have an equivalence of stabilities between good GPHs on $X$ and Higgs pairs on $Y$.

\begin{prop}\label{equivalence of stabilities between good GPHs and Higgs pairs}
A good GPH $(E,\phi_{E},F_{1}(E))$ is stable if and only if $(F,\phi_{F})$ is stable.
\end{prop}
\begin{proof}
Assume that $(E,\phi_{E},F_{1}(E))$ is stable of rank $n$ and degree $d$. Then
$$\deg F=\deg\nu_{*}E-\dim(\nu_{*}E\otimes\cc(p)/F_{1}(E))=\deg E+n-(2n-n)=\deg E$$
and $\phi_{E}=\nu^{*}\phi_{F}$. Let $K_{1}$ be a $\phi_{F}$-invariant subsheaf of $F$ of rank $r$. Let $K$ be the subbundle of $E$ generated by the image of $\nu^{*}K_{1}/\textrm{torsion}$ in $E$. Then $K$ is a $\phi_{E}$-invariant subbundle of $E$ of rank $r$ and $K_{1}$ is obtained by
$$0\to K_{1}\to\nu_{*}K\to\nu_{*}K\otimes\cc(p)/F_{1}(K)\to 0,$$
where $F_{1}(K)=F_{1}(E)\cap(K_{p_{1}}\oplus K_{p_{2}})$. Let $s=\dim F_{1}(K)$. Then we have
$$\deg K_{1}=\deg\nu_{*}K-\dim(\nu_{*}K\otimes\cc(p)/F_{1}(K))=\deg K+r-(2r-s)=\deg K+s-r.$$
Hence $\dss\frac{\deg K+s}{r}<\mu(E)+1$ if and only if $\dss\frac{\deg K_{1}}{r}<\mu(F)$.
\end{proof}

Combining Proposition \ref{equivalence of stabilities between GPBs and sheaves} and Lemma \ref{stability} with Proposition \ref{equivalence of semistabilities between good GPHs and Higgs pairs}, the semistability of good GPHs is simplified.

\begin{prop}\label{stability2}
A good GPH $(E,\phi_{E},F_{1}(E))$ is semistable if and only if $(E,F_{1}(E))$ is semistable.
\end{prop}

We complete the classification as follows.

\begin{thm}\label{coprime moduli and stable locus of noncoprime moduli of GPHs}
\begin{enumerate}
\item\label{coprime moduli and stable locus of noncoprime moduli of GPHs1} If $\gcd(n,d)=1$, then $\cM_{X}^{\GGPH}(n,d)\cong U_{X}^{\GPB}(n,d)\times\cc\cong\pp^{1}\times\cc$.
\item\label{coprime moduli and stable locus of noncoprime moduli of GPHs2} If $\gcd(n,d)>1$, then $\cM_{X}^{\GGPH}(n,d)^{s}$ is empty.
\end{enumerate}
\end{thm}
\begin{proof}
\begin{enumerate}
\item Note that if $(E,F_{1}(E))$ is stable, then any endomorphism of $(E,F_{1}(E))$ is a scalar (See Corollary 3.9 of \cite{Bh92}). The proof is similar to that of Theorem \ref{coprime moduli of Higgs pairs}. Proposition \ref{stability2} and the construction of $\cM_{X}^{\GGPH}(n,d)$ of \cite{Bh14} are applied to the proof.
\item The statement is an immediate consequence of Theorem \ref{noncoprime Higgs pairs} and Proposition \ref{equivalence of stabilities between good GPHs and Higgs pairs}.
\end{enumerate}
\end{proof}

Denote $\xymatrix{\pp^{1}\times\cc\ar[r]^{\cong\quad\quad}&\cM_{X}^{\GGPH}(n,d)}$ of Theorem \ref{coprime moduli and stable locus of noncoprime moduli of GPHs}-(\ref{coprime moduli and stable locus of noncoprime moduli of GPHs1}) by $\eta_{n,d}^{\GGPH}$.

Now we classify semistable Higgs pairs on $Y$. By Theorem \ref{noncoprime Higgs pairs} and using Jordan-H\"{o}lder filtration for Higgs pairs, we have the following observation.

\begin{lem}\label{splitting of ss Higgs pairs}
Assume $\gcd(n,d)=h$. Then any semistable Higgs pair $(F,\phi_{F})$ of rank $n$ and degree $d$ over $Y$ is S-equivalent to $(F_{1},\phi_{F_{1}})\oplus\cdots\oplus(F_{h},\phi_{F_{h}})$, where each $(F_{i},\phi_{F_{i}})$ is stable of rank $\dss\frac{n}{h}$ and degree $\dss\frac{d}{h}$.
\end{lem}

When $\gcd(n,d)=1$, there exists a universal family of stable Higgs pairs.

\begin{lem}\label{univ family of stable Higgs pairs}
Let $\gcd(n,d)=1$. There exists a universal family $\cE_{n,d}=(\cV_{n,d},\Phi_{n,d})$ of stable Higgs pairs of rank $n$ and degree $d$ parametrized by $Y\times\cc$.
\end{lem}
\begin{proof}
Consider the family of stable Higgs pairs $\cE_{n,d}=(\cV_{n,d},\Phi_{n,d})$ over $Y\times\cc$ such that
$$\cE_{n,d}|_{Y\times(y,t)}\cong((\cV_{n,d})_{y},\frac{t}{n}\otimes\id_{(\cV_{n,d})_{y}}).$$
Then for any $(y,t)\in Y\times\cc$, we have
$$\eta_{n,d}((y,t))=[\cE_{n,d}|_{Y\times(y,t)}]_{S},$$
where $[\cE_{n,d}|_{Y\times(y,t)}]_{S}$ is the S-equivalence class of $\cE_{n,d}|_{Y\times(y,t)}$.

Any family $\cF\to Y\times T$ induces canonically a morphism
$$\nu_{\cF}:T\to\cM_{Y}(n,d).$$
Then $f=(\eta_{n,d})^{-1}\circ\nu_{\cF}$ is a morphism $T\to Y\times\cc$ such that $\cF$ is S-equivalent to $f^{*}\cE_{n,d}$.
\end{proof}

If $\gcd(n,d)=h>1$, $\dss n'=\frac{n}{h}$ and $\dss d'=\frac{d}{h}$, then we can consider the family
$$\cE_{n,d}=\cE_{n',d'}\underset{Y}{\times}\overset{h}{\cdots}\underset{Y}{\times}\cE_{n',d'}$$
of polystable Higgs sheaves parametrized by $(Y\times\cc)\times\overset{h}{\cdots}\times(Y\times\cc)$.

\begin{rem}\label{permutation}
The action of the symmetric group $\fS_{h}$ on $Y\times\overset{h}{\cdots}\times Y$ induces an action of $\fS_{h}$ on $(Y\times\cc)\times\overset{h}{\cdots}\times(Y\times\cc)$. If $\om_{1}$ and $\om_{2}$ are two points of $(Y\times\cc)\times\overset{h}{\cdots}\times(Y\times\cc)$, the Higgs pairs $(\cE_{n,d})_{\om_{1}}$ and $(\cE_{n,d})_{\om_{2}}$ are S-equivalent if and only if $\om_{2}=\gamma\cdot\om_{1}$ for some $\gamma\in\fS_{h}$.
\end{rem}

The following map induced by $\cE_{n,d}$,
$$\nu_{\cE_{n,d}}:(Y\times\cc)\times\overset{h}{\cdots}\times(Y\times\cc)\to\cM_{Y}(n,d),$$
is surjective by Lemma \ref{splitting of ss Higgs pairs} and factors through $\Sym^{h}(Y\times\cc)$ by Remark \ref{permutation}. Now we complete the classification of semistable Higgs pairs as follows.

\begin{thm}\label{noncoprime moduli of Higgs pairs}
Let $\gcd(n,d)=h$.
\begin{enumerate}
\item\label{noncoprime moduli of Higgs pairs1} There exists a bijective morphism
$$\eta_{n,d}:\Sym^{h}(Y\times\cc)\to\cM_{Y}(n,d).$$
\item $\cM_{Y}(n,d)$ is irreducible.
\end{enumerate}
\end{thm}
\begin{proof}
\begin{enumerate}
\item $\nu_{\cE_{n,d}}$ induces a bijective morphism
$$\eta_{n,d}:\Sym^{h}(Y\times\cc)\to\cM_{Y}(n,d).$$
\item Since $\nu_{\cE_{n,d}}$ is continuous and $(Y\times\cc)\times\overset{h}{\cdots}\times(Y\times\cc)$ is irreducible, we get the result.
\end{enumerate}
\end{proof}

Finally we classify semistable good GPHs on $X$. By Theorem \ref{coprime moduli and stable locus of noncoprime moduli of GPHs}-(\ref{coprime moduli and stable locus of noncoprime moduli of GPHs2}) and using Jordan-H\"{o}lder filtration of GPHs, we have the following observation.

\begin{lem}\label{splitting of ss good GPHs}
Assume $\gcd(n,d)=h$. Then any semistable good GPH $(E,\phi_{E},F_{1}(E))$ of rank $n$ and degree $d$ over $X$ is S-equivalent to $(E_{1},\phi_{E_{1}},F_{1}(E_{1}))\oplus\cdots\oplus(E_{h},\phi_{E_{h}},F_{1}(E_{h}))$, where each $(E_{i},\phi_{E_{i}},F_{1}(E_{i}))$ is stable of rank $\dss\frac{n}{h}$ and degree $\dss\frac{d}{h}$.
\end{lem}

When $\gcd(n,d)=1$, there exists a universal family of stable good GPHs.

\begin{lem}\label{univ family of stable GGPHs}
Let $\gcd(n,d)=1$. There exists a universal family $\cE_{n,d}^{\GGPH}=(\cW_{n,d},\Phi_{n,d},F_{1}(\cW_{n,d}))$ of stable good GPHs of rank $n$ and degree $d$ parametrized by $X\times\cc$.
\end{lem}
\begin{proof}
Consider the family of stable good GPHs $\cE_{n,d}^{\GGPH}=(\cW_{n,d},\Phi_{n,d},F_{1}(\cW_{n,d}))$ over $X\times\cc$ such that
$$\cE_{n,d}^{\GGPH}|_{X\times(x,t)}\cong((\cW_{n,d})_{x},\frac{t}{n}\otimes\id_{(\cW_{n,d})_{x}},(F_{1}(\cW_{n,d}))_{x}).$$
Then for any $(x,t)\in X\times\cc$, we have
$$\eta_{n,d}^{\GGPH}((x,t))=[\cE_{n,d}^{\GGPH}|_{X\times(x,t)}]_{S},$$
where $[\cE_{n,d}^{\GGPH}|_{X\times(x,t)}]_{S}$ is the S-equivalence class of $\cE_{n,d}^{\GGPH}|_{X\times(x,t)}$.

Any family $\cF\to X\times T$ induces canonically a morphism
$$\nu_{\cF}:T\to\cM_{X}^{\GGPH}(n,d).$$
Then $f=(\eta_{n,d}^{\GGPH})^{-1}\circ\nu_{\cF}$ is a morphism $T\to X\times\cc$ such that $\cF$ is S-equivalent to $f^{*}\cE_{n,d}^{\GGPH}$.
\end{proof}

If $\gcd(n,d)=h>1$, $\dss n'=\frac{n}{h}$ and $\dss d'=\frac{d}{h}$, then we can consider the family
$$\cE_{n,d}^{\GGPH}=\cE_{n',d'}^{\GGPH}\underset{X}{\times}\overset{h}{\cdots}\underset{X}{\times}\cE_{n',d'}^{\GGPH}$$
of polystable good GPHs parametrized by $(\pp^{1}\times\cc)\times\overset{h}{\cdots}\times(\pp^{1}\times\cc)$.

The following map induced by $\cE_{n,d}^{\GGPH}$,
$$\nu_{\cE_{n,d}^{\GGPH}}:(\pp^{1}\times\cc)\times\overset{h}{\cdots}\times(\pp^{1}\times\cc)\to\cM_{X}^{\GGPH}(n,d),$$
is surjective by Lemma \ref{splitting of ss good GPHs} and factors through $\Sym^{h}(\pp^{1}\times\cc)$. We complete the classification of semistable good GPHs as follows.

\begin{thm}\label{noncoprime moduli of good GPHs}
Let $\gcd(n,d)=h$.
\begin{enumerate}
\item\label{noncoprime moduli of good GPHs1} There exists a bijective morphism
$$\eta_{n,d}^{\GGPH}:\Sym^{h}(\pp^{1}\times\cc)\to\cM_{X}^{\GGPH}(n,d).$$
\item $\Sym^{h}(\pp^{1}\times\cc)$ is the normalization of $\cM_{X}^{\GGPH}(n,d)$. 
\item $\cM_{X}^{\GGPH}(n,d)$ is irreducible.
\end{enumerate}
\end{thm}
\begin{proof}
\begin{enumerate}
\item $\nu_{\cE_{n,d}^{\GGPH}}$ induces a bijective morphism
$$\eta_{n,d}^{\GGPH}:\Sym^{h}(\pp^{1}\times\cc)\to\cM_{X}^{\GGPH}(n,d).$$
\item Since $\Sym^{h}(\pp^{1}\times\cc)$ is normal, the result follows by Zariski's main theorem.
\item Since $\nu_{\cE_{n,d}}^{\GGPH}$ is continuous and $(\pp^{1}\times\cc)\times\overset{h}{\cdots}\times(\pp^{1}\times\cc)$ is irreducible, we get the result.
\end{enumerate}
\end{proof}

\section{Hitchin map}\label{Hitchin map}
In this section, we give descriptions of all fibers of the Hitchin maps on $\cM_{X}^{\GGPH}(n,d)$ and $\cM_{Y}(n,d)$ respectively.

\begin{defn}[Section 5 of \cite{Bh14}]
The {\bf Hitchin map on $\cM_{X}^{\GGPH}(n,d)$} is defined by
$$H^{\GGPH}:\cM_{X}^{\GGPH}(n,d)\to A:=\bigoplus_{i=1}^{n}H^{0}(X,\cO_{X}),$$
$$(E,\phi_{E},F_{1}(E))\mapsto(a_{1}(\phi_{E}),\cdots,a_{n}(\phi_{E})),$$
where the characteristic polynomial $\det(\lambda-\phi_{E})$ of $(E,\phi_{E},F_{1}(E))$ is $\lambda^{n}+a_{1}(\phi_{E})\lambda^{n-1}+\cdots+a_{n}(\phi_{E})$.
\end{defn}

There is a result relating $\cM_{Y}(n,d)$ to $\cM_{X}^{\GGPH}(n,d)$.

\begin{prop}[Theorem 4.1 of \cite{Bh14}]\label{bir mor from GGPH to Higgs pair}
There exists a birational morphism
$$f:\cM_{X}^{\GGPH}(n,d)\to\cM_{Y}(n,d),(E,\phi_{E},F_{1}(E))\mapsto(F,\phi_{F}),$$
where $(F,\phi_{F})$ is given by (\ref{GGPH to Higgs pair}).
\end{prop}

\begin{rem}\label{surjectivity of f in coprime case}
If $\gcd(n,d)=1$, then $f$ is surjective (See Theorem 3 of \cite{Bh92}).
\end{rem}

Indeed $f$ is surjective in any case.

\begin{prop}\label{surjectivity of f}
$f$ is surjective.
\end{prop}
\begin{proof}
For each $(F,\phi_{F})\in\cM_{Y}(n,d)$, it follows from Lemma \ref{splitting of ss Higgs pairs} that
$$(F_{1},\phi_{F_{1}})\oplus\cdots\oplus(F_{h},\phi_{F_{h}}),$$
where each $(F_{i},\phi_{F_{i}})$ is stable of rank $\dss\frac{n}{h}$ and degree $\dss\frac{d}{h}$. By Remark \ref{surjectivity of f in coprime case}, there exists $(E_{i},\phi_{E_{i}},F_{1}(E_{i}))\in\cM_{X}^{\GGPH}(\frac{n}{h},\frac{d}{h})$ such that $f((E_{i},\phi_{E_{i}},F_{1}(E_{i})))=(F_{i},\phi_{F_{i}})$. Then we have the following commutative diagram of short exact sequences :
$$\xymatrix{0\ar[r]&\dss\bigoplus_{i=1}^{h}F_{i}\ar[r]\ar[dd]_{\dss\bigoplus_{i=1}^{h}\phi_{F_{i}}}&\nu_{*}\left(\dss\bigoplus_{i=1}^{h}E_{i}\right)\ar[r]\ar[dd]_{\dss\bigoplus_{i=1}^{h}\nu_{*}\phi_{E_{i}}}&\nu_{*}\left(\dss\bigoplus_{i=1}^{h}E_{i}\right)\otimes\cc(p)/\dss\bigoplus_{i=1}^{h}F_{1}(E_{i})\ar[r]\ar[dd]_{\left(\dss\bigoplus_{i=1}^{h}\nu_{*}\phi_{E_{i}}\right)_{p}}&0\\
&&&&\\
0\ar[r]&\dss\bigoplus_{i=1}^{h}F_{i}\ar[r]&\nu_{*}\left(\dss\bigoplus_{i=1}^{h}E_{i}\right)\ar[r]&\nu_{*}\left(\dss\bigoplus_{i=1}^{h}E_{i}\right)\otimes\cc(p)/\dss\bigoplus_{i=1}^{h}F_{1}(E_{i})\ar[r]&0.}$$
Hence
$$f((E_{1},\phi_{E_{1}},F_{1}(E_{1}))\oplus\cdots\oplus(E_{h},\phi_{E_{h}},F_{1}(E_{h})))=(F_{1},\phi_{F_{1}})\oplus\cdots\oplus(F_{h},\phi_{F_{h}}).$$
\end{proof}

\begin{prop}[Corollary 5.2 of \cite{Bh14}]\label{Hitchin maps are proper}
\begin{enumerate}
\item\label{Hitchin maps are proper1} $H^{\GGPH}:\cM_{X}^{\GGPH}(n,d)\to A$ is proper.
\item\label{Hitchin maps are proper2} $H^{\GGPH}$ defines a proper morphism
$$H:f(\cM_{X}^{\GGPH}(n,d))\to A,$$
where $f$ is given in Proposition \ref{bir mor from GGPH to Higgs pair}.
\end{enumerate}
\end{prop}

Proposition \ref{surjectivity of f} and Proposition \ref{Hitchin maps are proper}-(\ref{Hitchin maps are proper2}) imply the following statement.

\begin{cor}
$H^{\GGPH}$ defines a proper morphism
$$H:\cM_{Y}(n,d)\to A.$$
\end{cor}

This $H:\cM_{Y}(n,d)\to A$ is the {\bf Hitchin map on $\cM_{Y}(n,d)$}.

Now we describe all fibers of $H^{\GGPH}:\cM_{X}^{\GGPH}(n,d)\to A$ and $H:\cM_{Y}(n,d)\to A$. To describe all fibers of $H^{\GGPH}:\cM_{X}^{\GGPH}(n,d)\to A$, we follow the arguments of \cite[Section 5]{FGN14}.

Set $h=\gcd(n,d)$, $n'=\dss\frac{n}{h}$ and $d'=\dss\frac{d}{h}$. If $(E,\phi_{E},F_{1}(E))$ is a polystable good GPH of rank $n$ and degree $d$, then we have
$$(E,\phi_{E},F_{1}(E))\cong(E_{1},\phi_{E_{1}},F_{1}(E_{1}))\oplus\cdots\oplus(E_{h},\phi_{E_{h}},F_{1}(E_{h})),$$
where each $(E_{i},\phi_{E_{i}},F_{1}(E_{i}))$ is stable of rank $n'$ and degree $d'$ by Lemma \ref{splitting of ss good GPHs}. Then $(E_{i},\phi_{E_{i}})$ is stable by Proposition \ref{stability2} and $\phi_{E_{i}}=\dss\frac{t_{i}}{n'}\id_{E_{i}}$ where $\eta_{n',d'}^{\GGPH}((x_{i},t_{i}))=[(E_{i},\phi_{E_{i}})]_{S}$. Then
$$a_{1}(\phi_{E})=\sum_{i=1}^{h}t_{i},$$
$$\vdots$$
$$a_{n}(\phi_{E})=\left(\frac{t_{1}}{n'}\right)^{n'}\cdots\left(\frac{t_{h}}{n'}\right)^{n'}.$$
Denote the image of $H^{\GGPH}$ by $A_{n,d}$. If $D_{\bar{\lambda}}$ is the diagonal matrix with eigenvalues $\bar{\lambda}=(\lambda_{1},\cdots,\lambda_{h})$, then the following morphism
$$\alpha_{n,d}:\Sym^{h}\cc\to A_{n,d},$$
$$[\bar{t}]_{\fS_{h}}=[t_{1},\cdots,t_{h}]_{\fS_{h}}\mapsto(a_{1}(D_{\frac{1}{n'}\bar{t}}),\cdots,a_{n}(D_{\frac{1}{n'}\bar{t}}))$$
is bijective.

Define the projection
$$\pi_{h}:\Sym^{h}(\pp^{1}\times\cc)\to\Sym^{h}\cc,$$
$$[(x_{1},t_{1}),\cdots,(x_{h},t_{h})]_{\fS_{h}}\mapsto[t_{1},\cdots,t_{h}]_{\fS_{h}}.$$
Then the following diagram
$$\xymatrix{\Sym^{h}(\pp^{1}\times\cc)\ar[rr]^{\pi_{h}}\ar[d]_{\eta_{n,d}^{\GGPH}}^{\text{1:1}}&&\Sym^{h}\cc\ar[d]_{\text{1:1}}^{\alpha_{n,d}}\\
\cM_{X}^{\GGPH}(n,d)\ar[rr]^{H^{\GGPH}}&&A_{n,d}}$$
commutes.

The set $G=\{(t_{1},\cdots,t_{h})\in\cc^{h}\,|\,t_{i}\ne t_{j}\text{ if }i\ne j\}$ form an open dense subset of $\cc^{h}$. A point of $A_{n,d}$ is called {\bf generic} if it is the image under $\alpha_{n,d}$ of the $\fS_{h}$-orbit of some $\bar{t}_{g}\in G$. An arbirtrary point of $A_{n,d}$ is the image under $\alpha_{n,d}$ of a $h$-tuple of the form
$$\bar{t}_{a}=(t_{1},\overset{m_{1}}{\cdots},t_{1},\cdots,t_{l},\overset{m_{l}}{\cdots},t_{l}),$$
where $h=m_{1}+\cdots+m_{l}$.

\begin{prop}\label{fiber of pi h}
$$\pi_{h}^{-1}([\bar{t}_{g}]_{\fS_{h}})\cong\pp^{1}\times\overset{h}{\cdots}\times\pp^{1}$$
and
$$\pi_{h}^{-1}([\bar{t}_{a}]_{\fS_{h}})\cong\Sym^{m_{1}}\pp^{1}\times\overset{l}{\cdots}\times\Sym^{m_{l}}\pp^{1}\cong\pp^{m_{1}}\times\overset{h}{\cdots}\times\pp^{m_{l}}.$$
\end{prop}
\begin{proof}
The proof is same as that of Proposition 5.1 of \cite{FGN14}.
\end{proof}

\begin{cor}\label{fiber of h of GGPH}
The generic fiber of $H^{\GGPH}:\cM_{X}^{\GGPH}(n,d)\to A_{n,d}$ is set-theoretically isomorphic to $\pp^{1}\times\overset{h}{\cdots}\times\pp^{1}$. The fiber over an arbitrary point of the base is set-theoretically isomorphic to $\pp^{m_{1}}\times\overset{h}{\cdots}\times\pp^{m_{l}}$. The fiber over an arbitrary point of the base is isomorphic to $\pp^{1}$ for the case $\gcd(n,d)=1$.
\end{cor}
\begin{proof}
By Theorem \ref{noncoprime moduli of good GPHs}, there exists a bijective morphism from the fiber of $\pi_{h}$ to the fiber of $H^{\GGPH}$. So the first and the second statements follow from Proposition \ref{fiber of pi h}. The last statement follows from Theorem \ref{coprime moduli and stable locus of noncoprime moduli of GPHs}-(\ref{coprime moduli and stable locus of noncoprime moduli of GPHs1}).
\end{proof}

\begin{cor}\label{fiber of h of Higgs pair}
The generic fiber of $H:\cM_{Y}(n,d)\to A_{n,d}$ is set-theoretically isomorphic to $Y\times\overset{h}{\cdots}\times Y$. The fiber over an arbitrary point of the base is set-theoretically isomorphic to $\Sym^{m_{1}}Y\times\overset{l}{\cdots}\times\Sym^{m_{l}}Y$. The fiber over an arbitrary point of the base is isomorphic to $Y$ for the case $\gcd(n,d)=1$.
\end{cor}
\begin{proof}
Note that the normalization map $\nu:\pp^{1}\to Y$ induces the surjective map
$$g:\Sym^{h}(\pp^{1}\times\cc)\to\Sym^{h}(Y\times\cc),$$
$$[(x_{1},t_{1}),\cdots,(x_{h},t_{h})]_{\fS_{h}}\mapsto[(\nu(x_{1}),t_{1}),\cdots,(\nu(x_{h}),t_{h})]_{\fS_{h}}.$$
By Theorem \ref{noncoprime moduli of Higgs pairs}-(\ref{noncoprime moduli of Higgs pairs1}) and Theorem \ref{noncoprime moduli of good GPHs}-(\ref{noncoprime moduli of good GPHs1}), set-theoretically $f:\cM_{X}^{\GGPH}(n,d)\to\cM_{Y}(n,d)$ and $g$ coincide. The map $\pi_{h}$ and $g$ induce
$$\bar{\pi}_{h}:g(\Sym^{h}(\pp^{1}\times\cc))=\Sym^{h}(Y\times\cc)\to\Sym^{h}\cc,$$
$$[(y_{1},t_{1}),\cdots,(y_{h},t_{h})]_{\fS_{h}}\mapsto[t_{1},\cdots,t_{h}]_{\fS_{h}}.$$
Then we have the following commutative diagram
$$\xymatrix{\Sym^{h}(Y\times\cc)\ar[rr]^{\bar{\pi}_{h}}\ar[d]_{\eta_{n,d}}^{\text{1:1}}&&\Sym^{h}\cc\ar[d]_{\text{1:1}}^{\alpha_{n,d}}\\
\cM_{Y}(n,d)\ar[rr]^{H}&&A_{n,d}}$$

So $f(\pp^{1}\times\overset{h}{\cdots}\times\pp^{1})=Y\times\overset{h}{\cdots}\times Y$ and $f(\Sym^{m_{1}}\pp^{1}\times\overset{l}{\cdots}\times\Sym^{m_{l}}\pp^{1})=\Sym^{m_{1}}Y\times\overset{l}{\cdots}\times\Sym^{m_{l}}Y$. By Theorem \ref{noncoprime moduli of Higgs pairs} and Corollary \ref{fiber of h of GGPH}, there exists a bijective morphism from the fiber of $\bar{\pi}_{h}$ to the fiber of $H$. The last statement follows from Theorem \ref{coprime moduli of Higgs pairs}.
\end{proof}

\section{A flat degeneration}\label{degeneration}
Assume that $\gcd(n,d)=1$. In this section we consider a flat degeneration of the moduli space of stable Higgs bundles of rank $n$ and degree $d$ on an elliptic curve, which was constructed in \cite{BBN13,BBNarx13}.

We show that the moduli space of stable Higgs bundles of rank $n$ and degree $d$ on an elliptic curve degenerates to $Y\times\cc$ in a flat family. We also show that the fiber of the Hitchin map on the moduli space of stable Higgs bundles of rank $n$ and degree $d$ on an elliptic curve degenerates to $Y$ via the same flat family.

Let $Y^{(m)}$ be the curves which are semistably equivalent to $Y$, i.e. $X$ is a component of $Y^{(m)}$ and if $\nu:Y^{(m)}\to Y$ is the canonical morphism, $\nu^{-1}(p)$ is a chain $R$ of projective lines of length $m$, passing through $p_{1}$ and $p_{2}$.

\begin{defn}[\cite{NS99}]
A {\bf Gieseker vector bundle} of rank $n$ and degree $d$ on $Y^{(m)}$ is a vector bundle $E$ of rank $n$ and degree $d$ on $Y^{(m)}$ such that
\begin{itemize}
\item If $m\ge1$, then $E|_{R}$ is strictly standard,
\item $\nu_{*}E$ is a torsion-free sheaf on $Y$.
\end{itemize}
A Gieseker vector bundle $E$ on $Y^{(m)}$ is {\bf stable} if $\nu_{*}E$ is a stable torsion-free sheaf on $Y$.
\end{defn}

\begin{defn}[\cite{BBN13,BBNarx13}]
A {\bf Gieseker-Hitchin pair} of rank $n$ and degree $d$ on $Y^{(m)}$ is a pair $(E,\phi)$ such that $E$ is a Gieseker vector bundle of rank $n$ and degree $d$ on $Y^{(m)}$, $\phi$ is a global section of $\EEnd\,E$ and $\nu_{*}(E,\phi)$ is a Higgs pair on $Y$. A Gieseker-Hitchin pair $(E,\phi)$ on $Y^{(m)}$ is {\bf stable} if $\nu_{*}(E,\phi)$ is a stable Higgs pair on $Y$.
\end{defn}

Let $G_{Y}(n,d)$ be the moduli space of stable Gieseker vector bundles of rank $n$ and degree $d$ on $Y^{(m)}$ for some $0\le m\le n$ (See \cite[Theorem 1]{NS99}) and let $\cG_{Y}(n,d)$ the moduli space of stable Gieseker-Hitchin pairs of rank $n$ and degree $d$ on $Y^{(m)}$ for some $0\le m\le n$ (See \cite[Theorem 1]{BBN13} and \cite[Proposition 5.13]{BBNarx13}).

By Lemma \ref{stability}, the stability of Gieseker-Hitchin pairs can be simplified as follows.

\begin{lem}\label{stability of GH pairs}
A Gieseker-Hitchin pair $(E,\phi)$ on $Y^{(m)}$ is stable if and only if the underlying Gieseker vector bundle $E$ on $Y^{(m)}$ is stable.
\end{lem}

We can see that any endomorphism of a stable Gieseker-Hitchin pair is a scalar.

\begin{lem}\label{simpleness of endomorphism}
For a stable Gieseker vector bundle $E$ on $Y^{(m)}$, $\End\,E=\cc$.
\end{lem}
\begin{proof}
Since any finite dimensional division algebra over $\cc$ is $\cc$, it suffices to show that any nonzero $\varphi\in\End\,E$ is an isomorphism.

Consider a nonzero $\varphi\in\End\,E$. If $\varphi|_{X}:E|_{X}\to E|_{X}$ is zero, then $\varphi$ is zero at $p_{1}$ and $p_{2}$. Then $\varphi|_{R}:E|_{R}\to E|_{R}$ is zero, which is a contradiction. So $\varphi|_{X}:E|_{X}\to E|_{X}$ is nonzero. Since $\End(E|_{X})\cong\End(\nu_{*}E)$ by \cite[Remark 4]{NS99}, $\varphi|_{X}$ is an isomorphism. Then $\varphi$ is an isomorphism at $p_{1}$ and $p_{2}$. Since $(\det(E|_{R}))^{-1}\otimes(\det(E|_{R}))$ is trivial, $\det\varphi$ is nowhere zero and then $\varphi$ is an isomorphism. (See \cite[The proof of Proposition 3.1]{Gie84}).
\end{proof}

Now we define a flat family of Gieseker vector bundles, that of Gieseker-Hitchin pairs and their stabilities. Let $R$ be a discrete valuation ring with quotient field $K$ and residue field $\cc$. Let $T=\Spec R$, $\Spec K$ the generic point and $t_{0}$ the closed point of $T$. Let $Z\to T$ be a proper flat family such that the generic fiber $(Z_{K},{\bf0}:\Spec K\to Z_{K})$ is an elliptic curve and the closed fiber $Z_{t_0}\cong Y$. 

\begin{defn}[\cite{BBN13,BBNarx13}]
Let $Z^{(\md)}\to T$ is a flat morphism such that there exists a nonnegative integer $m$ satisfying that $(Z^{(\md)})_{t}=Y^{(m)}$ for each $t\in T$ with the commutative diagram
$$\xymatrix{Z^{(\md)}\ar[rr]^{\nu}\ar[rd]_{p_{T}}&&Z\ar[ld]\\
&T&}$$
where $\nu$ restricts to the morphism which contracts the chain $R$ of $\pp^{1}$'s on $Y^{(m)}$.
\begin{enumerate}
\item A {\bf Gieseker vector bundle on $Z^{(\md)}$} is a vector bundle $\cE_{T}$ on $Z^{(\md)}$ such that its restriction to the fiber $(Z^{(\md)})_{t}$ over $t\in T$ is a Gieseker vector bundle on $Y^{(m)}$ for some $m$.
\item A Gieseker vector bundle $\cE_{T}$ on $Z^{(\md)}$ is {\bf stable} if $\nu_{*}\cE_{T}$ is a family of stable torsion-free sheaves on $Z\to T$.
\item A {\bf Gieseker-Hitchin pair on $Z^{(\md)}$} is a pair $(\cE_{T},\varphi_{T})$ on $Z^{(\md)}$ such that $\cE_{T}$ is a vector bundle on $Z^{(\md)}$, $\varphi_{T}$ is a global section of $(p_{T})_{*}\EEnd\,\cE_{T}$ and its restriction to the fiber $(Z^{(\md)})_{t}$ over $t\in T$ is a Gieseker-Hitchin pair on $Y^{(m)}$ for some $m$.
\item A Gieseker-Hitchin pair $(\cE_{T},\varphi_{T})$ on $Z^{(\md)}$ is {\bf stable} if $\nu_{*}(\cE_{T},\varphi_{T})$ is a family of stable Higgs pairs on $Z\to T$.
\end{enumerate}
\end{defn}

Let $G_{Z/T}(n,d)\to T$ be the relative moduli space of stable Gieseker vector bundles of rank $n$ and degree $d$ on $Z^{(\md)}$ (See \cite[Theorem 2]{NS99}) and let $\cG_{Z/T}(n,d)\to T$ the relative moduli space of stable Gieseker-Hitchin pairs of rank $n$ and degree $d$ on $Z^{(\md)}$ (See \cite[Theorem 1]{BBN13} and \cite[Proposition 5.13]{BBNarx13}). Note that $G_{Z/T}(n,d)\to T$ and $\cG_{Z/T}(n,d)\to T$ are flat over $T$ (See \cite[Theorem 1]{BBN13}, \cite[Proposition 5.15]{BBNarx13} and \cite[Theorem 2]{NS99}). Moreover $G_{Z/T}(n,d)_{t}$ is the moduli space of vector bundles on the elliptic curve $Z_{t}$ for all $t\ne t_{0}$, $G_{Z/T}(n,d)_{t_{0}}\cong G_{Y}(n,d)$, $\cG_{Z/T}(n,d)_{t}$ is the moduli space of Higgs bundles on the elliptic curve $Z_{t}$ for all $t\ne t_{0}$ and $\cG_{Z/T}(n,d)_{t_{0}}\cong\cG_{Y}(n,d)$.

Now we describe $G_{Y}(n,d)$ and $\cG_{Y}(n,d)$ explicitly.

\begin{lem}\label{coprime moduli of GV}
$G_{Y}(n,d)\cong Y$.
\end{lem}
\begin{proof}
Let $G=G_{Z/T}(n,d)$. Then $G\to T$ is flat over $T$ such that $G_{t}\cong Z_{t}$ for all $t\ne t_{0}$ (See \cite{At57}) and $G_{t_{0}}\cong G_{Y}(n,d)$.

By \cite[Theorem 1 and Theorem 2]{NS99}, $G_{t_{0}}\cong G_{Y}(n,d)$ is a singular curve of arithmetic genus one.

If $\tG_{t_{0}}$ is the normalization of $G_{t_{0}}$ with the normalization morphism $\pi:\tG_{t_{0}}\to G_{t_{0}}$, then it follows from \cite[Corollary V.3.7 and Proposition V.3.8]{Ha77} that
$$p_{a}(\tG_{t_{0}})=p_{a}(G_{t_{0}})-\sum_{i=1}^{r}\frac{e_{i}(e_{i}-1)}{2},$$
where $p_{a}$ denotes the arithmetic genus and $e_{i}$ are the multiplicities of the infinitesimally near points of the singular points of $G_{t_{0}}$. Since $p_{a}(G_{t_{0}})=1$, then $p_{a}(\tG_{t_{0}})=0$, $r=1$ and $e_{i}=2$, which implies that $G_{t_{0}}$ has a unique ordinary double point, that is, a node or a cusp. Since $\tG_{t_{0}}$ is nonsingular and $p_{a}(\tG_{t_{0}})=0$, $\tG_{t_{0}}\cong\pp^{1}$.

By \cite[Lemma 3.1]{BD21}, the natural stratification of $G_{t_{0}}\cong G_{Y}(n,d)$
$$(G_{t_{0}})^{0}=G_{t_{0}}\supset(G_{t_{0}})^{1}\supset\cdots\supset(G_{t_{0}})^{n}\supset(G_{t_{0}})^{n+1}=\emptyset,$$
where $(G_{t_{0}})^{r+1}$ is the singular locus of $(G_{t_{0}})^{r}$ for every $0\le r\le n$ has the following property:
\begin{itemize}
\item $(G_{t_{0}})^{i}=\{x\in(G_{t_{0}})^{0}\,|\,\text{cardinality of the set }\pi^{-1}(x)\ge i+1\}\text{ for every }0\le i\le n$.
\item $(G_{t_{0}})^{i+1}$ is a Zariski-closed subvariety of $(G_{t_{0}})^{i}$ of pure codimension $1$, if non-empty.
\end{itemize}
Since $G_{t_{0}}$ is a curve, $(G_{t_{0}})^{2}=\emptyset$ and
$$(G_{t_{0}})^{1}=\{x\in(G_{t_{0}})^{0}\,|\,\text{cardinality of the set }\pi^{-1}(x)=2\}.$$
Hence the unique singular point of $G_{t_{0}}$ is a node and then $G_{t_{0}}\cong Y$.
\end{proof}

\begin{lem}\label{coprime moduli of GH pairs}
$\cG_{Y}(n,d)\cong G_{Y}(n,d)\times\cc\cong Y\times\cc$.
\end{lem}
\begin{proof}
The proof is same as that of Theorem \ref{coprime moduli of Higgs pairs} and Theorem \ref{flat degeneration}-(\ref{flat degeneration1}). We use Lemma \ref{stability of GH pairs}, Lemma \ref{simpleness of endomorphism} and Lemma \ref{coprime moduli of GV}.
\end{proof}

The Hitchin map on $\cG_{Z/T}(n,d)$ is defined as follows.

\begin{defn}[Definition 6.2 of \cite{BBNarx13}]
The {\bf Hitchin map on $\cG_{Z/T}(n,d)$} is defined by
$$H^{\GH}:\cG_{Z/T}(n,d)\to A_{T},$$
$$(E_{T},\varphi_{T})\mapsto(a_{1}(\varphi_{T}),\cdots,a_{n}(\varphi_{T})),$$
where $A_{T}\to T$ is the affine $T$-scheme representing the functor $S\mapsto\dss\bigoplus_{i=1}^{n}H^{0}(Z\times_{T}S,\cO_{Z\times_{T}S})$ and the characteristic polynomial $\det(\lambda-\varphi_{T})$ of $(E_{T},\varphi_{T})$ is $\lambda^{n}+a_{1}(\varphi_{T})\lambda^{n-1}+\cdots+a_{n}(\varphi_{T})$.
\end{defn}

\begin{prop}[Theorem 6.6 of \cite{BBNarx13}]\label{Hitchin map on flat family is proper}
$H^{\GH}:\cG_{Z/T}(n,d)\to A_{T}$ is proper over $T$.
\end{prop}

By \cite[Theorem 1 and Theorem 2]{NS99}, there exist proper and birational morphisms
$$\nu_{*}^{\GV}:G_{Y}(n,d)\to U_{Y}(n,d),E\mapsto\nu_{*}E$$
and
$$\nu_{*}^{\GV}:G_{Z/T}(n,d)\to U_{Z/T}(n,d),\cE_{T}\mapsto\nu_{*}\cE_{T},$$
where $U_{Z/T}(n,d)$ denotes the relative moduli space of stable torsion-free sheaves of rank $n$ and degree $d$ on $Z$.
By \cite[Corollary 5.14]{BBNarx13}, there exist proper and birational morphisms
$$\nu_{*}^{\GH}:\cG_{Y}(n,d)\to\cM_{Y}(n,d),(E,\phi)\mapsto\nu_{*}(E,\phi)$$
and
$$\nu_{*}^{\GH}:\cG_{Z/T}(n,d)\to\cM_{Z/T}(n,d),(\cE_{T},\varphi_{T})\mapsto\nu_{*}(\cE_{T},\varphi_{T})$$
where $\cM_{Z/T}(n,d)$ denotes the relative moduli space of stable Higgs pairs of rank $n$ and degree $d$ on $Z$.

Indeed we have the following observation.

\begin{prop}\label{pushforward is an id}
\begin{enumerate}
\item\label{pushforward is an id1} $\nu_{*}^{\GV}:G_{Y}(n,d)\to U_{Y}(n,d)$ is identified with the identity map $\id_{Y}:Y\to Y$.
\item\label{pushforward is an id2} $\nu_{*}^{\GH}:\cG_{Y}(n,d)\to\cM_{Y}(n,d)$ is identified with the identity map $\id_{Y\times\cc}:Y\times\cc\to Y\times\cc$.
\end{enumerate}
\end{prop}
\begin{proof}
By \cite[Proposition 3.2]{BD21} and Lemma \ref{coprime moduli of GV}, the singular locus of $G_{Y}(n,d)$ consists of stable Gieseker vector bundles on $Y^{(1)}$ and it corresponds to the single node $p$ of $Y$. Moreover the image of the singular locus of $G_{Y}(n,d)$ under $\nu_{*}^{GV}$ is exactly the singular locus of $U_{Y}(n,d)$ by \cite[Remark 2.3]{BB08}. Thus $\nu_{*}^{\GV}:G_{Y}(n,d)\to U_{Y}(n,d)$ is identified with the identity map $\id_{Y}:Y\to Y$. Moreover $\nu_{*}^{\GH}:\cG_{Y}(n,d)\to\cM_{Y}(n,d)$ is identified with the identity map $\id_{Y\times\cc}:Y\times\cc\to Y\times\cc$ by Theorem \ref{coprime moduli of Higgs pairs} and Lemma \ref{coprime moduli of GH pairs}.
\end{proof}

Lemma \ref{coprime moduli of GV} is relativized as follows.

\begin{lem}\label{relative coprime moduli of GV}
$G_{Z/T}(n,d)\cong Z$ as $T$-schemes.
\end{lem}
\begin{proof}
We first see that $U_{Z/T}(1,d)\cong Z$ as $T$-schemes. Let $T'=T\setminus\{t_{0}\}$. For any $T'$-scheme $S$, we have the following isomorphism
$$Z|_{T'}(S)\to U_{Z/T}(1,d)|_{T'}(S),\s\mapsto\cO_{Z|_{T'}}(\s)\otimes\cO_{Z|_{T'}}({\bf0})^{d-1}$$
on $S$-valued points, where ${\bf0}:S\to Z|_{T'}$ is a section which makes $(Z|_{T'},{\bf0})$ a flat family of elliptic curves. Then we have an isomorphism $\gamma:Z|_{T'}\to U_{Z/T}(1,d)|_{T'}$ by \cite[Lemma 5.7]{ACG11}. By \cite[Corollary 5.4]{ACG11}, $\gamma$ extends uniquely to an isomorphism $Z\to U_{Z/T}(1,d)$ over $T$.

Next we see that $U_{Z/T}(n,d)\cong Z$ as $T$-schemes. Consider the relative determinant morphism $\det_{Z/T}:U_{Z/T}(n,d)\to U_{Z/T}(1,d)$. Since $(\det_{Z/T})_{t}:U_{Z/T}(n,d)_{t}\to U_{Z/T}(1,d)_{t}$ is an isomorphism for all $t\in T'$, $\det_{Z/T}|_{T'}:U_{Z/T}(n,d)|_{T'}\to U_{Z/T}(1,d)|_{T'}$ is also an isomorphism by \cite[Lemma 5.7]{ACG11}. By \cite[Corollary 5.4]{ACG11}, $\det_{Z/T}|_{T'}$ extends uniquely to an isomorphism $\det_{Z/T}:U_{Z/T}(n,d)\to U_{Z/T}(1,d)$. Since we have seen that $U_{Z/T}(1,d)\cong Z$, we get an isomorphism $U_{Z/T}(n,d)\to Z$ over $T$.

Now we claim that $G_{Z/T}(n,d)\cong Z$ as $T$-schemes. Consider the morphism
$$\nu_{*}^{\GV}:G_{Z/T}(n,d)\to U_{Z/T}(n,d),\cE_{T}\mapsto\nu_{*}\cE_{T}.$$
For all $t\ne t_{0}$, $(\nu_{*}^{\GV})_{t}:G_{Z/T}(n,d)_{t}\to U_{Z/T}(n,d)_{t}$ is the identity map $\id_{Z_{t}}:Z_{t}\to Z_{t}$. By Proposition \ref{pushforward is an id}-(\ref{pushforward is an id1}), $(\nu_{*}^{\GV})_{t_{0}}:G_{Z/T}(n,d)_{t_{0}}\to U_{Z/T}(n,d)_{t_{0}}$ is also the identity map $\id_{Y}:Y\to Y$. Thus $\nu_{*}^{\GV}:G_{Z/T}(n,d)\to U_{Z/T}(n,d)$ is the identity map $\id_{Z}:Z\to Z$ by \cite[Corollary 5.4 and Lemma 5.7]{ACG11}. Since we have seen that $U_{Z/T}(n,d)\cong Z$ as $T$-schemes, we get $G_{Z/T}(n,d)\cong Z$ as $T$-schemes.
\end{proof}

We have a conclusion as follows.

\begin{thm}\label{flat degeneration}
\begin{enumerate}
\item\label{flat degeneration1} $\cG_{Z/T}(n,d)\cong G_{Z/T}(n,d)\times\cc\cong Z\times\cc$ as $T$-schemes.
\item $\nu_{*}^{\GH}:\cG_{Z/T}(n,d)\to\cM_{Z/T}(n,d)$ is identified with the identity map $\id_{Z\times\cc}:Z\times\cc\to Z\times\cc$.
\item The fiber of the Hitchin map $H^{\GH}$ on $\cG_{Z/T}(n,d)$ is isomorphic to $Z$.
\end{enumerate}
\end{thm}
\begin{proof}
\begin{enumerate}
\item By using the construction of $\cG_{Z/T}(n,d)$ in \cite[Section 5]{BBNarx13}, we show that $\cG_{Z/T}(n,d)\cong G_{Z/T}(n,d)\times\cc$.

Note that $U_{Z/T}(n,d)\cong R_{T}^{s}/\!/\PGL(d)$ for some $T$-scheme $R_{T}$. The Gieseker functor $\cG_{R_{T}}$ is represented by a $\PGL(d)$-invariant open subscheme $\cY$ of the $T$-scheme $\Hilb^{P_{1}}(Z\times_{T}\Gr(d,n))$ for some Hilbert polynomial $P_{1}$ and $G_{Z/T}(n,d)\cong\cY^{s}/\!/\PGL(d)$. Let $\Delta_{\cY}\subset Z\times_{T}\cY\times_{T}\Gr(d,n)$ be the universal object defining the functor $\cG_{R_{T}}$. The embedding $\Delta_{\cY}\subset Z\times_{T}\cY\times_{T}\Gr(d,n)$ gives the natural projection of $T$-schemes:
$$\xymatrix{\Delta_{\cY}\ar[rr]^{q}\ar[rd]_{p}&&\cY\ar[ld]^{r}\\
&T&}$$
Let $\cU$ be the universal vector bundle on $\Delta_{\cY}|_{\cY^{s}}$ obtained from the tautological quotient bundle on $\Gr(d,n)$. Applying \cite[Lemma 3.5]{Nit91} with $\cF=\EEnd\,\cU$, we get a linear scheme $\cY^{\GH}\to\cY^{s}$ given by
$$\cY^{\GH}=\bSpec\Sym_{\cO_{\cY^{s}}}(q_{*}\EEnd\,\cU)^{\vee}.$$
Since $q_{*}\EEnd\,\cU\cong\cO_{\cY^{s}}$ by Lemma \ref{simpleness of endomorphism} and the same argument of the proof of \cite[Lemma 4.6.3]{HL97}, we have $\cY^{\GH}\cong \cY^{s}\times\cc$. Lemma \ref{stability of GH pairs} implies that $(\cY^{\GH})^{s}\cong\cY^{s}\times\cc$.

Hence by the construction of \cite[Section 5]{BBNarx13} and Lemma \ref{relative coprime moduli of GV},
$$\cG_{Z/T}(n,d)\cong(\cY^{\GH})^{s}/\!/\PGL(d)\cong(\cY^{s}/\!/\PGL(d))\times\cc\cong G_{Z/T}(n,d)\times\cc\cong Z\times\cc.$$
\item We first show that $\cM_{Z/T}(n,d)\cong U_{Z/T}(n,d)\times\cc$ by the same argument as in the proof of Item (\ref{flat degeneration1}). Note that $U_{Z/T}(n,d)\cong R_{T}^{s}/\!/\PGL(d)$ for some $T$-scheme $R_{T}$ as mentioned. Let $\cU$ be the restriction of the universal sheaf to $Z\times_{T}R_{T}^{s}$. \cite[Lemma 3.5]{Nit91} with $\cF=\EEnd\,\cU$, the fact
$\pi_{R_{T}^{s}*}\EEnd\,\cU\cong\cO_{R_{T}^{s}}$ from \cite[Lemma 4.6.3]{HL97}, Lemma \ref{stability} and the proof of Lemma \ref{relative coprime moduli of GV} implies that
$$\cM_{Z/T}(n,d)\cong U_{Z/T}(n,d)\times\cc\cong Z\times\cc.$$
Now we describe $\nu_{*}^{\GH}:\cG_{Z/T}(n,d)\to\cM_{Z/T}(n,d)$. Since $\cG_{Z/T}(n,d)\cong G_{Z/T}(n,d)\times\cc$ in Item (\ref{flat degeneration1}) and $\cM_{Z/T}(n,d)\cong U_{Z/T}(n,d)\times\cc$ as shown previously, $\nu_{*}^{\GH}$ is exactly $\nu_{*}^{\GV}\times\id_{\cc}:G_{Z/T}(n,d)\times\cc\to U_{Z/T}(n,d)\times\cc$. In the proof of Lemma \ref{relative coprime moduli of GV}, we have seen that $\nu_{*}^{\GV}:G_{Z/T}(n,d)\to U_{Z/T}(n,d)$ is the identity map $\id_{Z}:Z\to Z$. Thus we get the statement.

\item By the definition of $H^{\GH}$ and $\cG_{Z/T}(n,d)\cong Z\times\cc$ in Item (\ref{flat degeneration1}), $H^{\GH}$ is identified with the projection $Z\times\cc\to\cc$ onto the second factor. Thus we get the result.
\end{enumerate}
\end{proof}

\end{document}